\documentclass[12pt]{article}

\usepackage{amsmath, epsfig, cite}
\usepackage{amssymb}
\usepackage{amsfonts}
\usepackage{latexsym}

\newtheorem{thm}{Theorem}[section]

\newtheorem{cor}[thm]{Corollary}
\newtheorem{conj}[thm]{Conjecture}
\newtheorem{lem}[thm]{Lemma}

\newcommand{\pf}{\noindent{\it Proof.} }

\def\N{{\mathbb N}}
\def\Z{{\mathbb Z}^+}

\numberwithin{equation}{section}

\newcommand{\qed}{{\hfill\rule{4pt}{7pt}}\medskip}

\begin{document}


\begin{center}
{\Large\bf New congruences for sums involving Ap\'ery numbers \\[5pt] or central Delannoy numbers}
\end{center}

\vskip 2mm \centerline{Victor J. W. Guo$^1$  and Jiang Zeng$^{2}$}
\begin{center}
{\footnotesize $^1$Department of Mathematics, East China Normal University,\\ Shanghai 200062,
 People's Republic of China\\
{\tt jwguo@math.ecnu.edu.cn,\quad http://math.ecnu.edu.cn/\textasciitilde{jwguo}}\\[10pt]
$^2$Universit\'e de Lyon; Universit\'e Lyon 1; Institut Camille
Jordan, UMR 5208 du CNRS;\\ 43, boulevard du 11 novembre 1918,
F-69622 Villeurbanne Cedex, France\\
{\tt zeng@math.univ-lyon1.fr,\quad
http://math.univ-lyon1.fr/\textasciitilde{zeng}} }
\end{center}


\vskip 0.7cm \noindent{\bf Abstract.}
The  Ap\'ery numbers $A_n$ and  central Delannoy numbers $D_n$ are defined by
$$A_n=\sum_{k=0}^{n}{n+k\choose 2k}^2{2k\choose k}^2, \quad
D_n=\sum_{k=0}^{n}{n+k\choose 2k}{2k\choose k}.
$$
Motivated by some recent work of Z.-W. Sun, we prove the  following  congruences:
\begin{align*}
\sum_{k=0}^{n-1}(2k+1)^{2r+1}A_k  &\equiv \sum_{k=0}^{n-1}\varepsilon^k (2k+1)^{2r+1}D_k \equiv 0\pmod n,
\end{align*}
where  $n\geqslant 1$, $r\geqslant 0$, and $\varepsilon=\pm1$. For $r=1$, we further show  that
\begin{align*}
\sum_{k=0}^{n-1}(2k+1)^{3}A_k &\equiv 0\pmod{n^3}, \quad \\
\sum_{k=0}^{p-1}(2k+1)^{3}A_k &\equiv p^3 \pmod{ 2p^6},
\end{align*}
where $p>3$ is a prime. The following congruence
\begin{align*}
\sum_{k=0}^{n-1} {n+k\choose k}^2{n-1\choose k}^2  \equiv 0 \pmod{n}
\end{align*}
plays an important role in our proof.


\vskip 2mm
\noindent{\it MR Subject Classifications}: 11A07, 11B65, 05A10

\section{Introduction}

The Ap\'ery numbers \cite{Apery} are given by
\begin{align}\label{eq:apery}
A_n
=\sum_{k=0}^{n}{n+k\choose 2k}^2{2k\choose k}^2.
\end{align}
Congruences for Ap\'ery numbers were studied by many people, see  Chowla et al. \cite{CCC}, Gessel \cite{Gessel}, and
Beukers \cite{Beukers}, for example.
The central Delannoy numbers (see \cite{CHV,Sulanke,Sun}) are defined by
\begin{align}\label{eq:delannoy}
D_n
=\sum_{k=0}^{n}{n+k\choose 2k}{2k\choose k}.
\end{align}
Recently, Sun \cite{Sun} has proved several remarkable congruences involving
Ap\'ery numbers or central Delannoy numbers.
In this paper we will prove some similar congruences related to Ap\'ery numbers and central Delannoy numbers,
and establish some new congruences on sums of binomial coefficients and $q$-binomial coefficients.
Let $\N$ denote the set of nonnegative integers and
$\Z$ the set of positive integers.
Our main results may be stated as follows.

\begin{thm}\label{thm:apery}
Let $r\in\N$ and $n\in\Z$. Then
\begin{align}
\sum_{k=0}^{n-1}(2k+1)k^r(k+1)^r A_k &\equiv 0 \pmod n, \label{eq:apery1}\\
\sum_{k=0}^{n-1}(2k+1)^{2r+1}A_k  &\equiv 0 \pmod n.   \label{eq:apery2}
\end{align}
\end{thm}

\begin{thm}\label{thm:del}
Let $r\in\N$ and $n\in\Z$. Then
\begin{align}
\sum_{k=0}^{n-1}\varepsilon^k (2k+1)k^r(k+1)^r D_k &\equiv 0 \pmod n, \label{eq:del1}\\
\sum_{k=0}^{n-1}\varepsilon^k (2k+1)^{2r+1}D_k  &\equiv 0 \pmod n, \label{eq:del2}
\end{align}
where $\varepsilon=\pm1$.
\end{thm}

When  $r=0$ the above theorems  were obtained by Sun~\cite{Sun}. For $r=1$, we have the following stronger result.
\begin{thm}\label{eq:aperyn3}
Let $n\in\Z$, and let $p>3$ be a prime. Then
\begin{align}
\sum_{k=0}^{n-1}(2k+1)^{3}A_k &\equiv 0\pmod{n^3}, \label{eq:akcubic}\\
\sum_{k=0}^{p-1}(2k+1)^{3}A_k &\equiv p^3 \pmod{ 2p^6}.  \label{eq:akcubic2}
\end{align}
\end{thm}

In order to prove Theorem~\ref{eq:aperyn3} we need to establish the congruence
\begin{align}
\sum_{k=0}^{n-1} {n+k\choose k}^2{n-1\choose k}^2  \equiv 0 \pmod{n}. \label{eq:particular1}
\end{align}
In fact, we have the following more general result.
\begin{thm}\label{thm:binom}
Let $a_1,\ldots,a_m,b_1,\ldots,b_m\in\N$ and $n\in\Z$. Then
\begin{align}
\sum_{k=0}^{n-1} {n-1\choose k}^2\prod_{i=1}^m {a_i+k\choose b_i+k}
&\equiv 0 \pmod{\gcd(a_1,\ldots,a_m,b_1,\ldots,b_m,n)}. \label{eq:noqversion}
\end{align}
\end{thm}

We shall prove Theorems 1.1--1.4 in the next four sections, respectively, along with some conjectures for further studies.

\section{Proof of Theorem \ref{thm:apery}}
 Let $(x)_0=1$ and $(x)_n=x(x+1)\cdots(x+n-1)$ for all  $n\in\Z$. We first  establish  three lemmas.

\begin{lem}\label{lem:coeff}
For all $k,m,r\in\N$, there exist $a_0(k,r), a_1(k,r),\ldots,a_r(k,r)\in\mathbb{Z}$  such that
\begin{align}
m^r(m+1)^r{m+k\choose 2k}
=\sum_{j=0}^{r}a_j(k,r) {m+k+j\choose 2k+2j} (2k+1)_{2j}.  \label{eq:coeff}
\end{align}
\end{lem}

\pf Given $k,r\in\N$, it is easy to see that
there exist integers $a_0(k,r), a_1(k,r),\ldots,a_r(k,r)$, independent of $x$, such that
\begin{align}\label{eq:newton}
x^r=\sum_{j=0}^r a_j(k,r) \prod_{i=1}^j (x-(k+i-1)(k+i)).
\end{align}
Substituting $x=m(m+1)$ in \eqref{eq:newton}, we get
\begin{align*}
m^r(m+1)^r
=\sum_{j=0}^r a_j(k,r) (m+k+1)_j (m-k-j+1)_j.
\end{align*}
Multiplying both sides by ${m+k\choose 2k}$ and noticing that
$$
{m+k\choose 2k} (m+k+1)_j (m-k-j+1)_j ={m+k+j\choose 2k+2j} (2k+1)_{2j},
$$
we obtain \eqref{eq:coeff}.  \qed
\begin{lem}\label{lem:gsun}
For all $a,k, n\in\N$, there holds
\begin{align}
\sum_{m=k}^{n-1}(2m+1){m+k\choose 2k}{m+k+a\choose 2k+2a}
=\frac{(n-k)(n-k-a)}{2k+a+1}{n+k\choose 2k}{n+k+a\choose 2k+2a}.\label{eq:gsun-00}
\end{align}
\end{lem}
Lemma \ref{lem:gsun} can be proved easily by induction on $n$ (see  Sun \cite[Lemma 2.1]{Sun} for the a=0 case).
We can also evaluate the left-hand side of \eqref{eq:gsun-00} automatically by Maple.

\begin{lem}\label{lem:congruence}
For all $a, k\in\N$ and $n\in\Z$, there holds
\begin{align*}
\frac{(n-k)(n-k-a)}{n(2k+a+1)}{n+k\choose 2k}{n+k+a\choose 2k+2a}{2k\choose k} (2k+1)_{2a} \in \N.
\end{align*}
\end{lem}

\pf For $a=0$, we have
\begin{align*}
\frac{(n-k)^2}{n(2k+1)}{n+k\choose 2k}^2{2k\choose k}
&={n-1\choose k}{n+k\choose k}{n+k\choose 2k+1}\in\N.
\end{align*}
For $a\geqslant  1$, we have
\begin{align*}
&\hskip -3mm
\frac{(n-k)(n-k-a)}{n(2k+a+1)}{n+k\choose 2k}{2k\choose k} (2k+1)_{2a} \\
&=(n-k-a){n+k\choose k}{n-1\choose k} (2k+1)_{a} (2k+a+2)_{a-1}\in\N.
\end{align*}
This completes the proof.    \qed

\noindent{\it Proof of Theorem {\rm\ref{thm:apery}.}} Substituting \eqref{eq:apery} for $A_m$ and
applying Lemmas \ref{lem:coeff} and \ref{lem:gsun}, we obtain
\begin{align*}
&\hskip -3mm\sum_{m=0}^{n-1}(2m+1)m^r(m+1)^r A_m \\
&=\sum_{j=0}^{r}\sum_{k=0}^{n-1} a_j(k,r){2k\choose k}^2(2k+1)_{2j}
\sum_{m=k}^{n-1} (2m+1){m+k\choose 2k}{m+k+j\choose 2k+2j} \\
&=\sum_{j=0}^{r}\sum_{k=0}^{n-1}a_j(k,r) {2k\choose k}C_j(k,n),
\end{align*}
where $a_0(k,r), a_1(k,r),\ldots,a_r(k,r)$  are integers determined by \eqref{eq:coeff} and
$$
C_j(k,n)={2k\choose k}(2k+1)_{2j}  \frac{(n-k)(n-k-j)}{2k+j+1}{n+k\choose 2k}{n+k+j\choose 2k+2j}.
$$
By Lemma \ref{lem:congruence},  we have $C_j(k,n)\equiv 0\pmod n$
and hence  \eqref{eq:apery1} holds. Since \begin{align} (2k+1)^{2r}
=(4k^2+4k+1)^r =\sum_{i=0}^r{r\choose i}4^{i}k^i(k+1)^i,
\label{eq:2k+1}
\end{align}
we deduce the congruence \eqref{eq:apery2} from \eqref{eq:apery1} immediately.  \qed

\section{Proof of Theorem \ref{thm:del}}
Similarly to the proof of Lemma \ref{lem:gsun}, we can check the following result (see \cite{Sun}).

\begin{lem}\label{lem:sun}
For all $k, n\in\N$, there hold
\begin{align*}
\sum_{m=0}^{n-1}(2m+1){m+k\choose 2k} &=\frac{(n-k)n}{k+1}{n+k\choose 2k}, \\
\sum_{m=0}^{n-1}(-1)^m (2m+1) {m+k\choose 2k}&=(-1)^{n-1}(n-k){n+k\choose 2k}.
\end{align*}
\end{lem}

\noindent{\it Proof of Theorem {\rm\ref{thm:del}.}}
Substituting \eqref{eq:delannoy}  for $D_m$ and
applying  Lemmas \ref{lem:coeff} and \ref{lem:sun}, we obtain
\begin{align*}
&\sum_{m=0}^{n-1}(2m+1)m^r(m+1)^r D_m\\
&=\sum_{j=0}^{r}\sum_{k=0}^{n-1}a_j(k,r) {2k\choose k}(2k+1)_{2j}  \frac{(n-k-j)n}{k+j+1}{n+k+j\choose 2k+2j} \\
&=n\sum_{j=0}^{r}\sum_{k=0}^{n-1}a_j(k,r) \frac{(n-k-j)}{k+j+1}{2k+2j\choose k+j} (k+1)_j^2 {n+k+j\choose 2k+2j},
\end{align*}
where $a_0(k,r),\ldots,a_r(k,r)$ are integers determined by \eqref{eq:coeff}.
Since $\frac{1}{k+j+1}{2k+2j\choose k+j}$ is an integer (a Catalan number),
we derive \eqref{eq:del1} for $\varepsilon=1$.
Similarly, we have
\begin{align*}
&\sum_{m=0}^{n-1}(-1)^m (2m+1)m^r(m+1)^r D_m\\
&=\sum_{j=0}^{r}\sum_{k=0}^{n-1} a_j(k,r){2k\choose k}(2k+1)_{2j} (-1)^{n-1} (n-k-j){n+k+j\choose 2k+2j} \\
&=\sum_{j=0}^{r}\sum_{k=0}^{n-1} a_j(k,r)(-1)^{n-1} (n-k-j){2k+2j\choose k+j} (k+1)_j^2 {n+k+j\choose 2k+2j}.
\end{align*}
Writing  $ (n-k-j){2k+2j\choose k+j} {n+k+j\choose 2k+2j}=n{n+k+j\choose n}{n-1\choose k+j}$, we
obtain  \eqref{eq:del1} for $\varepsilon=-1$.
Thus, the congruence \eqref{eq:del1} holds for $\varepsilon=\pm1$. By the relation \eqref{eq:2k+1},
we deduce the congruence \eqref{eq:del2} from \eqref{eq:del1} immediately.
\qed

The congruence \eqref{eq:del2} has the following refinements.
\begin{thm}\label{thm:delref}
Let $r\in\N$ and $n\in\Z$. Then
\begin{align*}
\sum_{k=0}^{n-1} (2k+1)^{2r+1}D_k
&\equiv n \pmod{2n},  \\
\sum_{k=0}^{n-1} (-1)^k (2k+1)^{2r+1}D_k
&\equiv
\begin{cases}
n,&\text{if $n$ is odd}\\[5pt]
0,&\text{otherwise}
\end{cases} \pmod{2n}.
\end{align*}
\end{thm}
\pf By the relation \eqref{eq:2k+1} and the congruence \eqref{eq:del1}, it is enough to prove the case $r=0$.
Substituting \eqref{eq:delannoy} for $D_k$ and  exchanging the order of summations,
we obtain, by applying Lemma \ref{lem:sun},  that
\begin{align*}
\sum_{k=0}^{n-1} (2k+1)D_k & = n\sum_{k=0}^{n-1}{n+k\choose n}{n\choose k+1}, \\
\sum_{k=0}^{n-1} (-1)^k (2k+1)D_k & = (-1)^n n\sum_{k=0}^{n-1}{n+k\choose n}{n-1\choose k},
\end{align*}
as already mentioned by Sun \cite{Sun}.
Noticing that
\begin{align}
\sum_{k=0}^{n-1}{n+k\choose n}{n\choose k+1}
&\equiv \sum_{k=0}^{n-1} (-1)^k {n+k\choose n}{n\choose k+1}
=(-1)^{n-1}  \pmod 2,  \nonumber \\
\sum_{k=0}^{n-1}{n+k\choose n}{n-1\choose k}
&\equiv \sum_{k=0}^{n-1}(-1)^k{n+k\choose n}{n-1\choose k}
=(-1)^{n-1} n  \pmod 2,  \label{eq:mod2}
\end{align}
we have
\begin{align*}
\sum_{k=0}^{n-1} (2k+1)D_k
&\equiv n \pmod{2n},   \\
\sum_{k=0}^{n-1} (-1)^k (2k+1)D_k
&\equiv
\begin{cases}
n,&\text{if $n$ is odd}\\[5pt]
0,&\text{otherwise}
\end{cases} \pmod{2n}.
\end{align*}
This completes the proof.   \qed

\begin{conj}Let $n$ be a power of $2$. Then
\begin{align*}
\sum_{k=0}^{n-1}(-1)^k(2k+1)^{3} D_k  &\equiv 2n^2 \pmod {n^3}.
\end{align*}
\end{conj}

\section{Proof of Theorem \ref{eq:aperyn3}}
We first prove two lemmas. The first one can be proved easily by induction on $n$, or be verified
automatically in Maple.
\begin{lem}\label{lem:2mk3}
For $0\leqslant  k\leqslant  n-1$ we have
\begin{align*}
\sum_{m=0}^{n-1}(2m+1)^3{m+k\choose 2k}^2=\frac{(n-k)^2(2n^2-k-1)}{k+1}{n+k\choose 2k}^2.
\end{align*}
\end{lem}
\begin{lem} For any prime $p>3$, we have
\begin{align}
\sum_{k=0}^{p-1}{p+k\choose k}^2{p-1\choose k}^2   &\equiv p \pmod{2p^4} , \label{eq:pksquare} \\
\sum_{k=0}^{p-1}{p+k\choose k+1}{p+k\choose k}{p-1 \choose k}^2 &\equiv 1 \pmod{2p^3}. \label{eq:pkpk+1}
\end{align}
\end{lem}
\pf  For $0\leqslant  k\leqslant  p-1$, we have
\begin{align*}
{p+k\choose k}^2{p-1\choose k}^2
&=\prod_{j=1}^k\frac{(p+j)^2(p-j)^2}{j^4}
=\prod_{j=1}^k\left(1-\frac{p^2}{j^2}\right)^2  \\
&\equiv 1-2p^2\sum_{j=1}^k \frac{1}{j^2} \pmod{p^4}.
\end{align*}
Therefore,
\begin{align}
\sum_{k=0}^{p-1}{p+k\choose k}^2{p-1\choose k}^2
&\equiv \sum_{k=0}^{p-1}\left(1-2p^2\sum_{j=1}^k\frac{1}{j^2}\right)  \nonumber\\
&= p-2p^2\sum_{j=1}^{p-1} \frac{p-j}{j^2} \pmod {p^4}.  \label{eq:ppkkpp-0}
\end{align}
By Wolstenholme's theorem (see \cite{HW}),  i.e., for $p>3$,
\begin{align}
\sum_{k=1}^{p-1}\frac{1}{k}\equiv 0 \pmod {p^2}, 
\qquad\sum_{k=1}^{p-1}\frac{1}{k^2}\equiv 0 \pmod {p},  \label{eq:psquare}
\end{align}
we obtain
\begin{align}
\sum_{k=0}^{p-1}{p+k\choose k}^2{p-1\choose k}^2   \equiv p \pmod{p^4}.  \label{eq:pkp4mod1}
\end{align}
Combining \eqref{eq:pkp4mod1} and \eqref{eq:mod2} with $n=p$, we obtain the congruence \eqref{eq:pksquare}.

Since ${p+k\choose k+1}=\frac{p}{k+1}{p+k\choose k}$, by \eqref{eq:ppkkpp-0} we have
\begin{align*}
{p+k\choose k+1}{p+k\choose k}{p-1 \choose k}^2
&\equiv \begin{cases}
1-2p^2\sum_{j=1}^k \frac{1}{j^2} , &\text{if $k=p-1$} \\
\frac{p}{k+1}\left(1-2p^2\sum_{j=1}^k \frac{1}{j^2} \right) &\text{if $0\leqslant  k\leqslant  p-2$}
\end{cases}  \pmod{p^4},
\end{align*}
and hence by \eqref{eq:psquare} we get
\begin{align*}
{p+k\choose k+1}{p+k\choose k}{p-1 \choose k}^2
&\equiv \begin{cases}
1, &\text{if $k=p-1$} \\
\frac{p}{k+1} &\text{if $0\leqslant  k\leqslant  p-2$}
\end{cases}  \pmod{p^3}.
\end{align*}
It follows that, for $p>3$,
\begin{align}
\sum_{k=0}^{p-1}{p+k\choose k+1}{p+k\choose k}{p-1 \choose k}^2
\equiv 1+\sum_{k=0}^{p-2}\frac{p}{k+1} \equiv 1  \pmod {p^3}  \label{eq:pk+1modp3}
\end{align}
by \eqref{eq:psquare}. On the other hand, we have
\begin{align}
\sum_{k=0}^{p-1}{p+k\choose k+1}{p+k\choose k}{p-1 \choose k}^2
&\equiv \sum_{k=0}^{p-1}{p+k\choose k+1}{p+k\choose k}{p-1 \choose k}  \nonumber \\
&=\sum_{k=0}^{p-1}{p+k\choose k}^2{p\choose k+1}  \nonumber \\
&\equiv \sum_{k=0}^{p-1}(-1)^k {p+k\choose k}{p\choose k+1}=(-1)^{p-1} \pmod 2.   \label{eq:pk+1mod2}
\end{align}
Combining \eqref{eq:pk+1modp3} and \eqref{eq:pk+1mod2}, we obtain the congruence \eqref{eq:pkpk+1}.
\qed

\noindent{\it Proof of Theorem {\rm\ref{eq:aperyn3}.}}
Substituting \eqref{eq:apery} for $A_m$ and
exchanging the order of summations, by Lemma \ref{lem:2mk3}, we obtain
\begin{align}
\sum_{m=0}^{n-1}(2m+1)^3 A_{m}
&=\sum_{k=0}^{n-1}{2k\choose k}^2 \frac{(n-k)^2(2n^2-k-1)}{k+1}{n+k\choose 2k}^2  \nonumber  \\
&=\sum_{k=0}^{n-1}{n+k\choose k}^2 \frac{(n-k)^2(2n^2-k-1)}{k+1}{n\choose k}^2  \nonumber\\
&=2n^3\sum_{k=0}^{n-1}{n+k\choose k+1}{n+k\choose k}{n-1\choose k}^2
-n^2 \sum_{k=0}^{n-1}{n+k\choose k}^2{n-1\choose k}^2. \label{eq:sum-cubic}
\end{align}
Applying \eqref{eq:particular1} to \eqref{eq:sum-cubic}, we immediately get the congruence \eqref{eq:akcubic}.
Now,  assume that $n=p>3$ is a prime in \eqref{eq:sum-cubic}. Then the congruence \eqref{eq:akcubic2} follows from
 \eqref{eq:pksquare} and \eqref{eq:pkpk+1}. \qed

\section{Proof of Theorem~\ref{thm:binom}}
We shall prove several $q$-versions of \eqref{eq:noqversion}. We first recall some $q$-notations and two fundamental  results.
The $q$-binomial coefficients and  $q$-integers are polynomials in $q$ with integer coefficients and are defined by
$$
{n\brack k}_{q}=\prod_{j=1}^{k}\frac{1-q^{n-j+1}}{1-q^{j}},\quad [n]_{q}=\frac{1-q^{n}}{1-q}.
$$
The following  two results are well known (see  \cite[(3.3.10)]{Andrews}  and \cite{Olive,De,GZ06}).
\begin{lem}[The $q$-Chu-Vandermonde formula] \label{prop:chuvan}
For nonnegative integers $m,n$ and $h$ there holds
\begin{align}
\sum_{k=0}^h {n\brack k}_{q}{m\brack h-k}_{q} q^{(n-k)(h-k)}={m+n\brack h}_q. \label{eq:q-chu}
\end{align}
\end{lem}

\begin{lem}[The $q$-Lucas Theorem]\label{prop:root}
Let $a, b, r, s$ and $n$ be  nonnegative  integers
such that $0\leqslant  b,s\leqslant  n-1$. Then, for any  $n$th  primitive root
of unity $\omega$, there holds
$$
{ad+b\brack rd+s}_{\omega}={a\choose r}{b\brack s}_\omega.
$$
\end{lem}

Given three polynomials $f(x)$, $g(x)$ and $q(x)$ in $\mathbb{Z}[x]$,
if $q(x)$ divides $f(x)-g(x)$ in $\mathbb{Z}[x]$ we write $f(x)\equiv g(x) \pmod{q(x)}$.
We shall prove the following $q$-versions of Theorem~\ref{thm:binom}.

\begin{thm}\label{thm:qbinom}
For any nonnegative integers $a_1,\ldots,a_m,b_1,\ldots,b_m$ and positive integer $n$, let
$d:=\gcd(a_1,\ldots,a_m,b_1,\ldots,b_m,n)$.
 Then
\begin{align}
\sum_{k=0}^{n-1} q^{k^2}{n-1\brack k}_q^2\prod_{i=1}^m {a_i+k\brack b_i+k}_q
&\equiv 0 \pmod{[d]_{q}},  \label{eq:qversion}  \\
 \sum_{k=0}^{n-1} q^{k^2+2k}{n-1\brack k}_q^2\prod_{i=1}^m {a_i+k\brack b_i+k}_q
&\equiv 0 \pmod{[d]_{q}},  \label{eq:qversion2} \\
  \sum_{k=0}^{n-1} q^{k}\prod_{i=1}^m {a_i+k\brack b_i+k}_q
&\equiv 0 \pmod{[d]_{q}},  \label{eq:qversion3} \\
  \sum_{k=0}^{n-1} q^{n-k-1}\prod_{i=1}^m {a_i+k\brack b_i+k}_q
&\equiv 0 \pmod{[d]_{q}},\label{eq:qversion4}
\end{align}
\end{thm}

\noindent{\it Proof.}
Let $d_1\in\Z$ such that $d_1|d$ and $d_1>1$, and let
$$
a_i=a_i'd_1,\ b_i=b_i'd_1,\ i=1,\ldots,m,\ \text{and}\ n=n'd_1.
$$
Suppose that $\omega$ is a primitive $d_1$-th root of unity.
Replacing $k$ by $jd_1+k$ with $0\leqslant  j\leqslant  n'-1$ and $0\leqslant k\leqslant d_1-1$, we obtain,
by the $q$-Lucas theorem,
\begin{align*}
\sum_{k=0}^{n-1} \omega^{k^2}{n-1\brack k}_\omega^2\prod_{i=1}^m {a_i+k\brack b_i+k}_\omega
&=\sum_{j=0}^{n'-1}\sum_{k=0}^{d_1-1} \omega^{(jd_1+k)^2}{n'd_1-1\brack jd_1+k}_\omega^2
\prod_{i=1}^m {a_i'd_1+jd_1+k\brack b_i'd_1+jd_1+k}_\omega  \\
&=\sum_{j=0}^{n'-1}\sum_{k=0}^{d_1-1} \omega^{k^2} {n'-1\choose j}^2 {d_1-1\brack  k}_\omega^2
\prod_{i=1}^m {a_i'+j\choose b_i'+j}  \\
&=\left(\sum_{j=0}^{n'-1}{n'-1\choose j}^2 \prod_{i=1}^m {a_i'+j\choose b_i'+j} \right)
\sum_{k=0}^{d_1-1} \omega^{k^2} {d_1-1\brack  k}_\omega^2.
\end{align*}
By applying the $q$-Chu-Vandermonde formula and the $q$-Lucas theorem, we have
$$
\sum_{k=0}^{d_1-1} \omega^{k^2} {d_1-1\brack  k}_\omega^2 ={2d_1-2\brack d_1-1}_{\omega}={d_1-2\brack d_1-1}_{\omega}=0.
$$
This proves that
\begin{align*}
\sum_{k=0}^{n-1} q^{k^2}{n-1\brack k}_q^2\prod_{i=1}^m {a_i+k\brack b_i+k}_q
 \equiv 0 \pmod{\Phi_{d_1}(q)}\quad\text{for any $d_1|d$ and $d_1>1$,}
\end{align*}
where $\Phi_n(q)$ is the $n$th cyclotomic polynomial in $q$.
Since
$$
\frac{1-q^d}{1-q}=\prod_{\substack{d_1|d\\ d_1>1}}\Phi_{d_1}(q),
$$
and the cyclotomic polynomials are pairwise relatively prime, we complete the proof of \eqref{eq:qversion}.
Similarly,  observing that
\begin{align*}
\sum_{k=0}^{d_1-1} \omega^{k^2+2k} {d_1-1\brack _k}_\omega^2
&=\omega^{-1}\sum_{k=0}^{d_1-1} \omega^{(d_1-k-1)^2} {d_1-1\brack _k}_\omega^2
=\omega^{-1}{2d_1-2\brack d_1-1}_{\omega}=0,   \\
\sum_{k=0}^{d_1-1} \omega^{k}  &= \sum_{k=0}^{d_1-1} \omega^{-k-1} =0,
\end{align*}
we can prove \eqref{eq:qversion2}--\eqref{eq:qversion4}.
\qed

To derive more consequences of Theorem \ref{thm:qbinom}, we need the following lemma.
\begin{lem}\label{lem:pol}
Let $P(x)\in\mathbb{Z}[x]$ and $d\in\Z$. If $P(n)\equiv  0\pmod d$ for all $n\in\N$,
then $P(m)\equiv  0\pmod d$ for all $m\in\mathbb{Z}$.
\end{lem}
\pf Just notice that $P(m+kd)\equiv P(m)\pmod d$ for all $m,k\in\mathbb{Z}$.  \qed

We now give the following generalization of Theorem~\ref{thm:binom}.
\begin{thm}\label{thm:gen-bino}
Let $a_1,\ldots,a_m\in\mathbb{Z}$, $b_1,\ldots,b_m\in\N$ and $n\in\Z$. Then
\begin{align*}
\sum_{k=0}^{n-1} (-1)^{mk} \prod_{i=1}^m {a_i-1\choose b_i+k}    \equiv 0 \pmod{\gcd(a_1,\ldots,a_m,b_1,\ldots,b_m,n)}.
\end{align*}
\end{thm}
\pf
Suppose that $d$ is a factor of $n$. Letting $q=1$ in  \eqref{eq:qversion3}, we obtain
\begin{align}
\sum_{k=0}^{n-1} \prod_{i=1}^m {x_id+k\choose y_id+k}    \equiv 0 \pmod{d}
\label{eq:aibi}
\end{align}
for all $x_1,\ldots,x_m,y_1,\ldots,y_m\in\N$.
Since $P(x_1):=\sum_{k=0}^{n-1} \prod_{i=1}^m {x_id+k\choose y_id+k}$
is a polynomial in $x_1$ with rational coefficients, there is a  positive integer $\alpha$ such that
$\alpha P(x_1)\in\mathbb{Z}[x_1]$. Then
\eqref{eq:aibi} is equivalent to
\begin{align}
\alpha P(x_1)   \equiv 0 \pmod{\alpha d}. \label{eq:pa1}
\end{align}
By Lemma \ref{lem:pol}, we see that \eqref{eq:pa1} is true for all $x_1\in\mathbb{Z}$, and so is \eqref{eq:aibi}.
By symmetry, we conclude that \eqref{eq:aibi} holds for all $x_1,\ldots,x_m\in\mathbb{Z}$.
Namely, for all  $a_1,\ldots,a_m\in\mathbb{Z}$, $b_1,\ldots,b_m\in\N$, there holds
\begin{align*}
\sum_{k=0}^{n-1} \prod_{i=1}^m {b_i-a_i+k\choose b_i+k}    \equiv 0 \pmod{\gcd(a_1,\ldots,a_m,b_1,\ldots,b_m,n)}.
\end{align*}
Noticing that ${b_i-a_i+k\choose b_i+k}=(-1)^{b_i+k}{a_i-1\choose b_i+k}$, we complete the proof.  \qed

Letting $b_i=0$ and $a_i=\pm n$ in Theorem \ref{thm:gen-bino}, we obtain
\begin{cor}Let $r,s,n\in\N$. Then
\begin{align*}
\sum_{k=0}^{n-1} {n+k\choose k}^r{n-1\choose k}^{2s}  &\equiv 0 \pmod{n},  \\
\sum_{k=0}^{n-1} (-1)^k {n+k\choose k}^r{n-1\choose k}^{2s+1}  &\equiv 0 \pmod{n},
\end{align*}
\end{cor}
In particular, we have
\begin{align}
\sum_{k=0}^{n}{n\choose k}^{2s}  \equiv 0 \pmod{(n+1)}.  \label{eq:rel-calkin}
\end{align}
It is worth mentioning that Calkin \cite[Proposition 3]{Calkin} has proved that
\begin{align}
\sum_{k=0}^{n}{n\choose k}^{2s}  \equiv 0 \pmod{p}  \label{eq:calkin}
\end{align}
if $p$ is a prime such that $\frac{n}{m}<p<\frac{n+1}{m}+\frac{n+1-m}{m(2ms-1)}$ for some $m\in\Z$.
It is clear that \eqref{eq:rel-calkin} and \eqref{eq:calkin} are different and can not
be deduced from each other.

Letting $b_i=0$ and $a_i=-n$ or $a_i=-2n$ in Theorem \ref{thm:gen-bino}, we obtain
\begin{cor}Let $r,s,n\in\N$. Then
\begin{align*}
\sum_{k=0}^{n-1}(-1)^{(r+s)k} {n-1\choose k}^{r}{2n-1\choose k}^s &\equiv 0 \pmod{n}.
\end{align*}
\end{cor}
Note that Chamberland and Dilcher \cite{CD} have studied divisibility properties of the following similar sums:
$$
\sum_{k=0}^{n}\varepsilon^{k} {n\choose k}^{r}{2n\choose k}^s,
$$
where $\varepsilon=\pm 1$.

\begin{conj}Let $a_1,\ldots,a_m\in\mathbb{Z}$, $b_1,\ldots,b_m,r\in\N$ and $n\in\Z$. Then
\begin{align*}
\sum_{k=0}^{n-1} (-1)^{mk} \varepsilon^k (2k+1)k^r(k+1)^r\prod_{i=1}^m {a_i-1\choose b_i+k}
  &\equiv 0 \pmod{\gcd(a_1,\ldots,a_m,b_1,\ldots,b_m,n)},
\end{align*}
\end{conj}
where $\varepsilon=\pm 1$.

\begin{conj}Let $n,r\in\Z$. Then
\begin{align*}
\sum_{k=0}^{n-1}{n-1\choose k}^{2r}
&\equiv
\begin{cases}
n,&\text{if $n=2^a$}\\[5pt]
0,&\text{otherwise}
\end{cases} \pmod{2n}, \\
\sum_{k=0}^{n-1}{2n-1\choose k}^{2r}
&\equiv
\begin{cases}
n,&\text{if $n=2^a$}\\[5pt]
0,&\text{otherwise}
\end{cases} \pmod{2n}.
\end{align*}
\end{conj}

\begin{conj}Let $n,r\in\N$ and $s,t\in\Z$. Then
\begin{align*}
\sum_{k=0}^{n-1}(-1)^{kt}{n+k\choose k}^{s}{n-1\choose k}^{t}
&\equiv
\begin{cases}
0,&\text{if $2\mid n$ and $2\nmid (s+t)$}\\[5pt]
n,&\text{otherwise}
\end{cases} \pmod{2n}.
\end{align*}
\end{conj}

\begin{lem}Let $n\in\N$. Then
\begin{align}\label{eq:chu-van-bis}
\sum_{k=0}^n{n\brack k}_q^2 q^{k^2-k}=2{2n-1\brack n}_q.
\end{align}
\end{lem}

\pf  By the $q$-Chu-Vandermonde formula \eqref{eq:q-chu}, we have
\begin{align*}
\sum_{k=0}^n{n\brack k}_q {n-1\brack k}_q  q^{k^2}&={2n-1\brack n}_q, \\
\sum_{k=0}^n{n\brack k}_q  {n-1\brack k-1}_q q^{k(k-1)}&={2n-1\brack n}_q.
\end{align*}
The result then follows from the relation ${n\brack k}_q={n-1\brack k}_q q^k+{n-1\brack k-1}_q $. \qed

\begin{thm}Let $n\in\Z$. Then
\begin{align}
\sum_{k=0}^{n-1} q^{k^2-k}{n+k\brack k}_q^2{n-1\brack k}_q^2
\equiv
\begin{cases}
2I(n/2)\frac{1-q^n}{1-q^2}, &\text{if $n$ is even}\\[5pt]
0, &\text{otherwise}
\end{cases}\pmod{  [n]_q},  \label{eq:fin-qbino}
\end{align}
where $I(n)=\frac{1}{n}\sum_{k=0}^{n-1} {n+k\choose k}^2{n-1\choose k}^2$.
\end{thm}
\pf
Let $d_1|n$ ($d_1>1$) and $n'=n/d_1$. Suppose that $\omega$ is a primitive $d_1$-th root of unity.
Similarly to the proof of Theorem \ref{thm:qbinom}, we have
\begin{align*}
\sum_{k=0}^{n-1} \omega^{k^2-k}{n-1\brack k}_\omega^2{n+k\brack k}_\omega^2
&=\sum_{j=0}^{n'-1}{n'-1\choose j}^2 {n'+j\choose j}^2
\sum_{k=0}^{d_1-1} \omega^{k^2-k} {d_1-1\brack _k}_\omega^2. \\
&=2{2d_1-3\brack d_1-1}_\omega\sum_{j=0}^{n'-1}{n'-1\choose j}^2 {n'+j\choose j}^2
\quad\text{by}\ \eqref{eq:chu-van-bis}.
\end{align*}
If $n$ is odd, then $d_1\geqslant  3$, and by the $q$-Lucas theorem,
$$
{2d_1-3\brack d_1-1}_\omega={d_1-3\brack d_1-1}_\omega=0,
$$
which means that
\begin{align*}
\sum_{k=0}^{n-1} q^{k^2-k}{n+k\brack k}_q^2{n-1\brack k}_q^2
\equiv
0 \pmod{[n]_q}.
\end{align*}
If $n$ is even, then it is easy to check that (for $d_1=2$ or $d_1\geqslant  3$)
\begin{align*}
\sum_{k=0}^{n-1} \omega^{k^2-k}{n-1\brack k}_\omega^2{n+k\brack k}_\omega^2
-2I(n/2)\sum_{k=0}^{(n-2)/2}\omega^{2k}=0,
\end{align*}
which implies the first congruence in \eqref{eq:fin-qbino}.
This completes the proof.    \qed

\begin{conj}\label{conj:qsquare}
Let $n$ be any  power of a prime  $p$. Then
\begin{align*}
\sum_{k=0}^{n-1} q^{(n-k)^2}{n+k\brack k}_q^2{n-1\brack k}_q^2
\equiv q^{(n-1)^2}[n]_q\pmod{ [p]_{q^{n/p}}^2}.
\end{align*}
\end{conj}

It is easy to see that Conjecture \ref{conj:qsquare} is true for $q=1$
by the congruences  \eqref{eq:particular1} and \eqref{eq:pksquare}.

\vskip 5mm \noindent{\bf Acknowledgments.} This work was partially
supported by the Fundamental Research Funds for the Central
Universities,  Shanghai Rising-Star Program (\#09QA1401700) and the
National Science Foundation of China (\#10801054).


\begin{thebibliography}{99}
\small \setlength{\itemsep}{-.8mm}
\bibitem{Andrews} G.E. Andrews, The Theory of Partitions,
Cambridge University Press, Cambridge, 1998.

\bibitem{Apery}R. Ap\'ery, Irrationalit\'e de $\zeta(2)$ et $\zeta(3)$, Ast\'erisque 61 (1979), 11--13.

\bibitem{Beukers}F. Beukers, Another congruence for the Ap\'ery numbers, J. Number Theory 25
(1987), 201--210.


\bibitem{Calkin}N.J. Calkin, Factors of sums of powers of binomial coefficients,
Acta Arith. 86 (1998), 17--26.

\bibitem{CHV}J.S. Caughman, C.R. Haithcock and J.J.P. Veerman, A note on lattice chains
and Delannoy numbers, Discrete Math. 308 (2008), 2623--2628.

\bibitem{CD}M. Chamberland and K. Dilcher,
Divisibility properties of a class of binomial sums, J. Number Theory 120 (2006), 349--371.

\bibitem{CCC}S. Chowla, J. Cowles and M. Cowles, Congruence properties of Ap\'ery numbers,
J. Number Theory 12 (1980), 188-190.

\bibitem{De} J. D\'esarm\'enien,  Un analogue des congruences de Kummer pour les $q$-nombres d'Euler,
European J. Combin. {3} (1982), 19--28.


\bibitem{Gessel}I. Gessel, Some congruences for Ap\'ery Numbers, J. Number Theory 14 (1982), 362--368.


\bibitem{GZ06}V.J.W. Guo and J. Zeng, Some arithmetic properties of the $q$-Euler numbers and $q$-Sali\'e
numbers, European J. Combin. 27 (2006), 884--895.

\bibitem{HW}G.H. Hardy and E.M. Wright, An Introduction to the Theory of Numbers, 5th Editon, the Clarendon Press,
Oxford University Press, New York, 1979, p.~88.

\bibitem{Olive}G. Olive, Generalized powers, Amer. Math. Monthly {72} (1965), 619--627.

\bibitem{Sulanke}R.A. Sulanke, Objects counted by the central Delannoy numbers,
J. Integer Seq. 6 (2003),  Article 03.1.5.

\bibitem{Sun}Z.-W. Sun, On sums of Ap\'ery polynomials and related congruences, preprint, arXiv:1101.1946v3.


\bibitem{Sunopen}Z.-W. Sun, Open conjectures on congruences, preprint, arXiv:0911.5665v39.




\end{thebibliography}
\end{document}